\documentclass{amsart}

\usepackage{amsthm,amssymb}
\usepackage[pagebackref,colorlinks,linkcolor=red,citecolor=blue,urlcolor=blue,hypertexnames=true]{hyperref}
\usepackage[alphabetic]{amsrefs}

\newcommand{\cC}{\mathcal C}
\newcommand{\cD}{\mathcal D}
\newcommand{\N}{\mathbb N}
\newcommand{\R}{\mathbb R}
\newcommand{\Z}{\mathbb Z}
\newcommand{\res}{{\rm res}}
\newcommand{\ran}{{\rm ran}}
\newcommand{\dom}{{\rm dom}}

\theoremstyle{plain}
\newtheorem{theorem}{Theorem}[section]
\newtheorem{corollary}[theorem]{Corollary}
\newtheorem{lemma}[theorem]{Lemma}
\theoremstyle{definition}
\newtheorem{definition}[theorem]{Definition}

\theoremstyle{remark}
\newtheorem{remark}[theorem]{Remark}
\newtheorem{example}[theorem]{Example}

\begin{document}
\title[Bounded cohomology and L$^2$-invariants]{Low degree bounded cohomology and L$^2$-invariants for negatively curved groups}


\author{Andreas Thom}
\address{Andreas Thom, Mathematisches Institut,
Bunsenstr. 3-5\\ 37073 G\"ottingen, Germany.}
\email{thom@uni-math.gwdg.de}
\urladdr{http://www.uni-math.gwdg.de/thom/}

\begin{abstract}
We study the subgroup structure of discrete groups that 
share cohomological properties which resemble non-negative
curvature. Examples include all Gromov hyperbolic groups.

We provide strong restrictions on the possible $s$-normal subgroups of a 'negatively curved' group.
Another result says that the image of a group, which is boundedly generated by a finite set
of amenable subgroups, in a group, which admits a proper 
quasi-1-cocycle into the regular representation, has to be amenable. These results extend to a certain class of randomorphisms in the sense of Monod. 
\end{abstract}

\maketitle

\tableofcontents

\section{Introduction} \label{intro}

This note is a companion to a paper by J.\ Peterson and the author \cite{petthom}, to which we refer for background and notation. Cohomology and bounded cohomology of an infinite group with coefficients in the left regular representation have proved to be useful tools to understand properties of the group. For the study of $\ell^2$-homology and cohomology
we refer to the book by W.\ L\"uck \cite{lueck}, for information about bounded cohomology the standard reference
is the book by N.\ Monod \cite{MR1840942}.

Non-vanishing of the second bounded cohomology with coefficients in the left regular representation is the key condition in the work of Burger-Monod \cite{MR1911660}, Monod-Shalom \cites{MR2153026,ms} and Mineyev-Monod-Shalom \cite{MR2081428} on rigidity theory. In \cite{petthom},
Peterson and the author studied non-vanishing of the first cohomology with coefficients in the left regular representation
and derived results about the subgroup structure. Now, we link first $\ell^2$-cohomology and second bounded cohomology via an exact sequence of $LG$-modules, where $LG$ denotes the group von Neumann algebra. Moreover, we extend the methods of \cite{petthom} to apply to a wider class of groups including all
non-elementary Gromov hyperbolic groups. 
The key notion here is the notion of quasi-1-cocycle which appeared already
in several places including Monod's foundational work \cite{MR1840942}.
\medskip

The organization of the article is as follows. Section \ref{intro} is the Introduction. After studying quasi-1-cocycles and the first quasi-cohomology group  in connection with L\"uck's dimension theory in Section \ref{intro} and \ref{dimthe}, we introduce a class of groups $\cD_{\rm reg}$ which is closely related to the class $\cC_{\rm reg}$, which was studied in
\cites{MR2153026,ms}. Examples of groups in the class $\cD_{\rm reg}$ include non-elementary Gromov hyperbolic groups and all groups with a positive first $\ell^2$-Betti number. In the sequel we prove two main results. In Section \ref{snormal} we prove that an s-normal subgroup of a group in $\cD_{\rm reg}$ is
again in $\cD_{\rm reg}$ using methods from \cite{petthom}. In Section \ref{boundgen} we show that all group homomorphisms from a group, which is boundedly generated by a finite set of amenable subgroups to a group, which admits a proper quasi-1-cocycle into an infinite sum of the left regular representation, have amenable image. If one adds property (T) to the assumptions on the source group, the image has to be finite. Since the combination of bounded generation and property (T) applies to many lattices in higher rank Lie groups, this gives a new viewpoint towards some well-known results in the field.

In the last section, we study the class of groups which admit such proper quasi-1-cocycles and show that it is closed under formation of free products and a notion of $\ell^2$-orbit equivalence. The methods imply that the results of Section \ref{boundgen} extend to certain randomorphisms in the sense of Monod, see \cite{MR2275641}.

\vspace{0.2cm}

Throughout the article, $G$ will be a discrete countable group, $\ell^2 G$ denotes the Hilbert space with basis $G$ endowed with the left regular representation.

\begin{definition}
Let $\pi \colon G \to U(H_{\pi})$ be a unitary representation of $G$.
A map $c\colon G \to H_{\pi}$ is called a \emph{quasi-1-cocycle}, if and only if the map
$$G\times G \ni (g,h) \to \pi(g)c(h) - c(gh) + c(g) \in H_{\pi}$$
is uniformly bounded on $G \times G$. We denote the vectorspace of quasi-1-cocycles with values in $H_{\pi}$
by $QZ^1(G,H_{\pi})$. The subspace of uniformly bounded quasi-1-cocycles is denoted by $QB^1(G,H_{\pi})$.
\end{definition}

In analogy to the definition of $H^1(G,H_\pi)$, we define the following:

\begin{definition}
The first quasi-cohomology of $G$ with coefficients in the unitary $G$-representation $H_{\pi}$ is defined by
$$QH^1(G,H_{\pi}) = QZ^1(G,H_{\pi})/QB^1(G,H_{\pi}).$$
\end{definition}

The relevance of $QH^1(G,H_\pi)$ becomes obvious in the next theorem, which links low degree cohomology and
second bounded cohomology with coefficients in $H_{\pi}$.

\begin{theorem} \label{sequence}
Let $G$ be a discrete countable group and $H_\pi$ be a unitary $G$-representation. 
There exists an exact sequence as follows:
$$0 \to H^1(G,H_\pi) \to QH^1(G,H_\pi) \stackrel{\delta}{\to} H^2_b(G,H_\pi) \to H^2(G,H_\pi).$$
Here, $\delta$ denotes the Hochschild coboundary, which maps a 1-chain to a 2-cocycle.
\end{theorem}
\begin{proof}
The proof is contained in \cite{MR1840942}. However, since this sequence is relevant for our work, we give a short argument. Clearly, the first non-trivial map is injective, since the trivial 1-cocycles are precisely the bounded ones, by the
Bruhat-Tits fixed point lemma.
The kernel of $\delta$ is given precisely by those 1-chains which are cycles and define elements in $H^1(G,H_\pi)$.
It remains to prove exactness at $H^2_b(G,H_\pi)$. This is obvious, since the elements which are mapped
to zero are precisely those which are coboundaries of 1-chains. Thus, these are precisely those which are in the image of the coboundary map $\delta$.
\end{proof}

It is clear, that homological algebra can fit every comparison map between additive functors, which is defined on a suitable \emph{chain level}, into a long exact sequence. Hence, the above exact sequence can be extended to the right, see \cite{MR1840942}.
Since we will not need this extension we do not dwell on this.

\section{Second bounded cohomology and dimension theory} \label{dimthe}

The following theorem was first observed for finitely generated groups
by Burger-Monod in \cite{MR1911660} as a consequence
of the new approach to bounded cohomology developed by N.\ Monod in \cite{MR1840942}.
Later, the result was extended to arbitrary countable groups by V.\ Kaimanovich \cite{MR2006560}.
Note, that we do not state the most general form of the result, but rather a version which we
can readily apply to the problems we study.

\begin{theorem} \label{reph2} Let $G$ be a discrete countable group and $H_{\pi}$ be a separable dual Banach $G$-module. There exists a standard probability space $S$ with a Borel $G$-action, leaving the probability measure quasi-invariant, 
such that there exists a natural isometric isomorphism 
$$H^2_b(G,H_{\pi}) \cong ZL^{\infty}_{\rm alt}(S^3,H_{\pi})^G.$$
Here, $ZL^{\infty}_{\rm alt}(S^3,H_{\pi})$ denotes the space of alternating 3-cocycles 
with values in $H_{\pi}$.
\end{theorem}

The naturality of the isomorphism immediately implies that, in case $H_{\pi}$ carries an additional module structure over a ring $R$ (which commutes with the $G$-action), the isomorphism is an isomorphism of $R$-modules.
One highly non-trivial consequence of the preceding theorem 
(from our point of view) is that $H^2_b(G,H_{\pi})$ can be viewed as a subspace
of $L^{\infty}(S^3,H_{\pi})$, for some \emph{standard} probability space $S$. 
The following corollary makes use of this fact.

In the sequel, we are freely using L\"uck's dimension function for $LG$-modules, which is defined for \emph{all} 
modules over the group von Neumann algebra $LG$. For details about its definition see \cite{lueck}. Note that
$\ell^2 G$ carries a commuting right $LG$-module structure, which induces $LG$-module structures on all 
(quasi-)cohomological invariants of it.

\begin{definition}
A $LG$-module $M$ is called \emph{rank separated} if for every non-zero element $\xi \in M$
$$[\xi] \stackrel{\rm def}{=}1- \sup\{\tau(p)\mid p^2=p^*=p \in LG, \xi p =0\} > 0.$$
\end{definition}

For more information on the notion of rank, we refer to \cite{thomrank}. The only consequence we need is the following lemma.

\begin{lemma} 
A $LG$-module is rank separated if and only if every non-zero
$LG$-submodule has a positive dimension.
\end{lemma}

\begin{corollary} \label{bcsep}
Let $G$ be a discrete countable group and $K \subset G$ a subgroup.
Then, $H^2_b(K,\ell^2G)$ is a rank separated $LG$-module. In particular, every non-zero element generates a sub-module
of positive dimension.
\end{corollary}
\begin{proof}
If we apply Theorem \ref{reph2} to $K$, it follows that $H^2_b(K,\ell^2 G) \subset L^{\infty}(S^3,\ell^2G)$. 
Let $\xi \in H^2_b(K,\ell^2 G)$.
and assume that there exists a sequence 
of projections $p \in LG$, such that $p_n \uparrow 1$ and $\xi p_n=0$, for all $n \in \N$. 
If $\xi p_n =0$, then for a co-null set $X_n \subset S^3$, we have $\xi(x)p_n=0$, for all $x \in X_n$. Clearly,
the intersection $\cap_{n \in \N} X_n$ is still co-null and hence $\xi(x) =0$, for almost all $x\in S^3$. This implies $\xi =0$.
\end{proof}

\begin{corollary} \label{nonam}
Let $G$ be a discrete countable group and $K \subset G$ a subgroup. 
The $LG$-module $$QH^1(K,\ell^2G)$$ is rank separated if and only if
$K$ is non-amenable.
\end{corollary}
\begin{proof}
In view of the exact sequence
$$0 \to H^1(K,\ell^2 G) \to QH^1(K,\ell^2 G) \to H^2_b(K,\ell^2 G),$$
the $LG$-module $QH^1(K,\ell^2 G)$ is rank separated if both $H^1(K,\ell^2 G)$ and $H^2_b(K,\ell^2 G)$ are
rank separated. Hence, the proof is finished by Corollary \ref{bcsep} and Hulanicki's Theorem, see also
the proof of Corollary $2.4$ in \cite{petthom}.
\end{proof}

Following \cite{ms}, we denote  by $\cC_{\rm reg}$ the class of groups, for which $$H^2_b(G,\ell^2 G) \neq 0.$$
This class was studied extensively in \cite{ms}, and strong results about rigidity and superrigidity were obtained.
An a priori slightly different class is of importance in the results we obtain.

\begin{definition}
We denote by $\cD_{\rm reg}$ the class of groups with $$\dim_{LG} QH^1(G,\ell^2G) \neq 0.$$ 
\end{definition}

\begin{remark}
Neither of the possible inclusions between $\cC_{\rm reg}$ and $\cD_{\rm reg}$ is known to hold. Any positive result in this direction would be very interesting. Both of the inclusions seem to
be likely. Indeed, it is very likely that
$\dim_{LG}QH^1(G,\ell^2G)$ and $\dim_{LG} H^2_b(G,\ell^2 G)$ can only take the values $0$ or $\infty$. 
At least if $G$ is of type $FP_2$, this would imply $G \in \cC_{\rm reg} \Leftrightarrow G \in \cD_{\rm reg}$. 
Moreover, it would give the implication $$ \beta_1^{(2)}(G) \in (0,\infty) \Rightarrow G \in \cC_{\rm reg},$$
which seems natural. However, we did not succeed in proving the required restriction on the values of
the dimension.
\end{remark}

Both, $\cC_{\rm reg}$ and $\cD_{\rm reg}$ consist of groups which all remember some features of negatively curved metric spaces. It is therefore permissible to call these groups 'negatively curved'. We will see this more
directly in the examples below.

\begin{lemma} \label{lem2}
Let $G$ be a countable discrete group. If the group $G$ is in $\cD_{\rm reg}$, then 
either the first $\ell^2$-Betti number of $G$
or the second bounded cohomology of $G$ with coefficients in $\ell^2 G$ does not vanish.
The converse holds if $\beta_2^{(2)}(G)=0$.
\end{lemma}
\begin{proof}
This is immediate from Theorem \ref{sequence} and Corollary \ref{bcsep}.
\end{proof}

A large class of groups in $\cD_{\rm reg}$ is provided by a result from \cite[Thm.\ 3]{MR2081428}.

\begin{theorem}
All non-elementary hyperbolic groups are in $\cD_{\rm reg}$.
\end{theorem}
\begin{proof}
In \cite{MR2081428}, it was shown that all hyperbolic groups have non-vanishing $QH^1(G,\ell^2G)$. In fact, it was shown
that there exists some element in $QH^1(G,\ell^2 G)$ which even maps non-trivially to $H^2_b(G,\ell^2 G)$. 
It follows from Corollary \ref{nonam}, that all non-elementary (i.e.\ non-amenable) hyperbolic groups are in $\cD_{\rm reg}$.
\end{proof}

In \cite{MR2153026}, the class $\cC_{\rm reg}$ is studied more extensively. 
One result we want to mention is the following:

\begin{theorem}[Corollary $7.6$ in \cite{MR2153026}]
Let $G$ be a discrete group acting non-elementarily 
and properly by isometries on some proper {\rm CAT(-1)} space.
Then, $$H^2_b(G,\ell^2 G) \neq 0.$$
\end{theorem}

Note that in view of Lemma \ref{lem2}, the preceding theorem provides examples of 
groups in $\cD_{\rm reg}$ as soon as the second $\ell^2$-Betti number of the corresponding group
vanishes.

\section{Non-existence of infinite $s$-normal subgroups} \label{snormal}

The following notion of normality was studied by Peterson and the author \cite{petthom} in connection with a non-vanishing first $\ell^2$-Betti number. The definition of $s$-normality goes back to the seminal work of S.\ Popa, who studied similar definitions in
\cite{popa}.

\begin{definition}
Let $G$ be a discrete countable group. 
An infinite subgroup $K \subset G$ is said to be $s$-normal, if $gKg^{-1} \cap K$ is infinite for all $g \in G$.
\end{definition}

\begin{example}
The inclusions $GL_n(\Z) \subset GL_n({\mathbb Q})$ and $\Z = \langle a \rangle \subset \langle a,b \mid ba^pb^{-1} = a^q \rangle =BS_{p,q}$ are inclusions of $s$-normal subgroups.
\end{example}

Given Banach space valued functions $f,g\colon X \to B$, defined on a set $X$, we write $f \asymp g$ if 
the function $$X \ni x \mapsto \|f(x)-g(x)\| \in \R$$ is uniformly bounded on $X$.

A unitary $G$-representation is said to be \emph{strongly mixing}, if $\langle g\xi,\eta \rangle \to 0$,
for $g \to \infty$. The following lemma is the key observation which leads to our first main results.

\begin{lemma} \label{mainlem1}
Let $G$ be a discrete countable group and let $K \subset G$ be an infinite $s$-normal subgroup.
Let $H_\pi$ be a strongly mixing unitary representation of the group $G$.
The restriction map
$$\res^G_K\colon QH^1(G,H_\pi) \to QH^1(K,H_\pi)$$
is injective.
\end{lemma}
\begin{proof}
Let $c\colon G \to H_\pi$ be a quasi-1-cocycle, which is bounded on $K$. Thus, $c(k) \asymp 0$ as a function
of $k\in K$. We compute
$$c(gkg^{-1}) \asymp (1-gkg^{-1})c(g) + gc(k) \asymp (1-gkg^{-1})c(g),$$
as a function of $(g,k) \in G \times K$. For $k \in g^{-1}Kg \cap K$, we obtain
$$\|(1-gkg^{-1})c(g)\| \leq \|c(gkg^{-1})\| + C'  \leq C,$$
for some constants $C',C \geq 0$.
Since $K \subset G$ is $s$-normal, the subgroup $g^{-1}Kg \cap K$ is infinite. Now, 
since $H_\pi$ is strongly mixing,
we conclude that $$2\|c(g)\|^2 = \lim_{k \to \infty} \|(1-gkg^{-1})c(g)\|^2 \leq C^2.$$ 
Hence, $g \to \|c(g)\|$ is uniformly bounded by $2^{-1/2} C$. This proves the claim.
\end{proof}

The following theorem is a non-trivial consequence about the subgroup structure for groups in the class $\cD_{\rm reg}$.
It is our first main result.

\begin{theorem}
Let $G$ be a discrete countable group in $\cD_{\rm reg}$
and let $K \subset G$ be an infinite
$s$-normal subgroup. The group $K$ satisfies at least one of the following properties:
\begin{enumerate}
\item[i] the first $\ell^2$-Betti number of $K$ does not vanish, or
\item[ii] the second bounded cohomology of $K$ with coefficients in $\ell^2K$ does not vanish.
\end{enumerate}
In particular, $K$ can neither be amenable nor a product of infinite groups.
\end{theorem}
\begin{proof}
By Lemma \ref{mainlem1}, the restriction map 
$$QH^1(G,\ell^2G) \to QH^1(K,\ell^2 G)$$
is injective. One easily sees that $QH^1(QK,\ell^2 K) \subset QH^1(K,\ell^2 G)$ generates a rank-dense $LG$-submodule in $QH^1(H,\ell^2 G)$. Hence $QH^1(K,\ell^2 K)$ cannot be zero-dimensional and we see that 
$K$ is in $\cD_{\rm reg}$ and the claim follows from
Lemma \ref{lem2}. We conclude that $K$ can neither be amenable or a product of infinite groups, since both classes of groups have vanishing first $\ell^2$-Betti number (see \cite{lueck}) and vanishing second bounded cohomology with coefficients in the left regular representation, see \cite{ms}.
\end{proof}

\begin{remark}
The result easily extends to $ws$-normal subgroups, see \cite{petthom}
\end{remark}

\section{Bounded generation and finiteness theorems} \label{boundgen}

It has been observed by many people that boundedly generated groups and non-elementary 
hyperbolic groups are opposite extremes in geometric group theory. In this section we support this view by showing that
there are essentially no group homomorphisms from a boundedly generated with property (T) to a Gromov hyperbolic group. Later, in Section \ref{obeq} we can even extend this result to a suitable class of randomorphisms in the sense of Monod, see \cite{MR2275641}.

A group $G$ is said to be \emph{boundedly generated} by a subset $X$, if
there exists $k \in \N$, such that each element of $G$ is a product of less than $k$ elements from $X \cup X^{-1}$. We say that $G$ is boundedly generated by a finite set of subgroups $\{G_i, i\in I\}$ if $G$ is 
boundedly generated by the set $\cup_{i \in I } G_i$.

\begin{lemma} \label{bounded}
Let $G$ be a non-amenable group which is boundedly generated by a finite set of amenable subgroups.
Then, the group $QH^1(G,\ell^2G^{\oplus \infty})$ is zero.
\end{lemma}
\begin{proof}
We view $\ell^2G^{\oplus \infty} \cong \ell^2(G \times \Z)$ and consider it as a $L(G\times \Z)$-module.
In view of Corollary \ref{nonam}, we can assume that $QH^1(G,\ell^2G^{\oplus \infty})$ is
rank separated. Hence, given an arbitrary 
element $$c \in QH^1(G,\ell^2 (G \times \Z)),$$ in order to show that it is zero,
we have to provide a sequence of projection $p_n \in L(G \times \Z)$, such that 
$p_n \uparrow 1$ and $cp_n=0$, for all $n \in \N$.

Let $G$ be boundedly generated by amenable subgroups $G_1,\dots,G_n$. The restriction of a quasi-1-cocycle
onto $G_i$ is almost bounded, i.e.\ there exists a projection $q_i \in L(G \times \Z)$ 
of trace $\tau(q_i) \geq 1-{\varepsilon}/n$, such that
$cq_i$ is bounded on $G_i$. Setting $p = \inf_{1 \leq i \leq n} q_i$, we obtain a projection $p$ with trace $\tau(p) \geq 1-
\varepsilon$, such that $cp$ is bounded on $G_i$, for all $1 \leq i \leq n$. The cocycle identity and bounded generation imply that $cp$ is bounded on the whole of $G$. Hence $cp=0 \in QH^1(G,\ell^2(G \times \Z))$. The sequence $p_n$ is  constructed by choosing $\varepsilon < 1/n$. This proves the claim.
\end{proof}
As we have seen, non-elementary Gromov hyperbolic groups are in $\cD_{\rm reg}$ but more is true:
\begin{lemma} \label{proper1}
Let $G$ be a Gromov hyperbolic group. There exists a proper quasi-1-cocycle on $G$ with values in
$\ell^2 G^{\oplus \infty}$.
\end{lemma}
\begin{proof} This is an immediate consequence of the proof of Thm.\ $7.13$ in \cite{MR2153026}, which
builds on Mineyev's work on equivariant bicombings, Thm.\ $10$ in \cite{MR1866802}.
\end{proof}

The following theorem is the second main result of this article. A cocycle version of it will
be presented in the last section as Corollary \ref{coromain2}.

\begin{theorem} \label{main2}
Let $G$ be a group which admit a proper quasi-1-cocycle into $\ell^2 G^{\oplus \infty}$ and let $H$ be a group which is boundedly generated by
a finite set of amenable subgroups. Then, every group homomorphism $\phi\colon H \to G$ has amenable image.
\end{theorem}
\begin{proof} We may assume that $\phi$ is injective, since any quotient of a group which is boundedly generated
by amenable groups is of the same kind. 
If the quotient is non-amenable, then the restriction of the proper quasi-1-cocycle coming from Lemma \ref{proper1}
has to be bounded on $H$, by Lemma \ref{bounded}. Indeed, $\ell^2H^{\oplus \infty} \cong \ell^2G^{\oplus \infty}$ 
as unitary $H$-representations using a coset decomposition and hence Lemma \ref{bounded} applies.
However, the quasi-1-cocycle is unbounded on any infinite subset. 
This is a contradiction, since $H$ follows to be finite and hence amenable. 
\end{proof}
\begin{remark}
The result applies in particular in the case when $G$ is Gromov hyperbolic. In this case we can even conclude that
the image is finite or virtually cyclic, since all amenable subgroups of a Gromov hyperbolic group are finite or virtually cyclic.
\end{remark}

From the above theorem we can derive the following corollary.

\begin{corollary} \label{coro2}
Let $G$ be a group which admits a proper quasi-1-cocycle into $\ell^2 G^{\oplus \infty}$ 
and let $H$ be a group which
\begin{enumerate}
\item[i] is boundedly generated by a finite set of amenable subgroups, and
\item[ii] has property (T) of Kazhdan-Margulis.
\end{enumerate}
Then, every group homomorphism $\phi\colon H \to G$ has finite image.
\end{corollary}
\begin{proof}
Any quotient of a property (T) group has also property (T). However, the image of $\phi$ is amenable
by Theorem \ref{main2} and the only amenable groups with
property (T) are finite. This finishes the proof.
\end{proof}

\begin{remark} Note that results like the preceding Corollary are well-known if one assumes the target to
be a-T-menable, whereas here: many Gromov hyperbolic groups have property (T). 
Examples of groups $H$ which satisfy assumptions $i$ and $ii$ in Corollary \ref{coro2}
include $SL_n(\Z)$ for $n \geq 3$ and many other lattices in higher rank 
semi-simple Lie groups, see \cite{MR1044049}. Conjecturally, all irreducible, non-cocompact 
lattices in higher rank Lie groups share these properties. In \cite{MR1911660}, it was shown that
higher rank lattices in certain algebraic groups over local fields have property (TT) of Monod, see Theorem $13.4.1$
in \cite{MR1840942} and the definitions therein. A similar proof can be carried out in this situation.

Note that the groups which satisfy the conditions $i$ and $ii$ do not always satisfy property (TT) of Monod.
Indeed, in an appendix of \cite{MR2197372} Monod-R\'{e}my construct boundedly 
generated groups (in fact lattices in higher rank semi-simple Lie groups) with property (T) which fail 
to have property (QFA) of Manning (see \cite{MR2197372}) and property (TT) of Monod. 
Hence, the plain quasification of property (T), 
which would yield some version of property (TT) of Monod (and should of course also imply property (QFA) by an extension of Watatani's proof, see \cite{wat}) is too strong to hold for all lattices in higher 
rank semi-simple Lie group. Hence, the combination of conditions $i$ and $ii$ is maybe the appropriate set of 
conditions that encodes the way in which higher rank lattices satisfy a \emph{strong form} of property (T).
\end{remark}

\begin{remark}
The mechanism of properness vs.\ boundedness works in the context of ordinary first cohomology with coefficients
in $\ell^p$-spaces as well. In fact, G.\ Yu \cite{yu} provided proper $\ell^p$-cocycles for hyperbolic groups, whereas
Bader-Furman-Gelander-Monod \cite{bfgm} studied the necessary strengthening of property (T) for higher rank lattices.
A combination of the Theorem B in \cite{bfgm} and Yu's result allows to conclude Corollary \ref{coro2} for lattices in certain algebraic groups (see the assumption of Theorem $B$ in \cite{bfgm}). However, even for this special case, our approach seems more elementary, just using the notion of quasi-1-cocycle.
\end{remark}

\section{Groups with proper quasi-1-cocycles}

\subsection{Subgroups and free products}

Every quasi-1-cocycle is close to one for which:
$$c(g^{-1}) = -g^{-1}c(g)$$
holds on the nose. Indeed, $\tilde{c}(g) = \frac12\left(c(g) - g c(g^{-1})\right)$ is anti-symmetric and only
bounded distance away from $c$. Note also that $\tilde{c}$ is proper if and only $c$ is proper. 
We call the quasi-1-cocycles which satisfy this additional property \emph{anti-symmetric}.

\begin{lemma} \label{freeprod}
Let $G,H$ be discrete countable groups and 
let $H_{\pi}$ be a unitary representation of $G \ast H$. Moreover, let $c_1\colon G \to H_{\pi}$ and
$c_2 \colon H \to H_{\pi}$ be anti-symmetric quasi-1-cocycles. Then, there is a natural anti-symmetric quasi-1-cocycle 
$$c = (c_1 \ast c_2) \colon G \ast H \to H_\pi$$
which extends $c_1$ and $c_2$.
\end{lemma}
\begin{proof}
Let $w=g_1h_1g_2h_2 \dots g_n h_n$ be a reduced element in $G \ast H$ (i.e. only $g_1$ or $h_n$ might be trivial).
We define
$$c(w) = c_1(g_1) + g_1 c_2(h_2) + g_1h_1c_1(g_2) + \dots + g_1h_1 \dots g_n c_2(h_n).$$
Clearly, $c$ is anti-symmetric, i.e. $c(w^{-1}) = - w^{-1}c(w)$ just by construction and using that $c_1$ and $c_2$ were
anti-symmetric. Let us now check that it is indeed a quasi-1-cocycle.
Let $w_1$ and $w_2$ be elements of $G \ast H$ and assume that $w_1 = w'_1 r$ and $w_2 = r^{-1} w'_2$, such that
the products $w_1'w_2', w_1'r$, and $r^{-1}w'_2$ are reduced, in the sense that the block length drops at most by one.

Then, the identities 
\begin{eqnarray*}
c(w_1'w_2') &=& w_1'c(w_2') + c(w_1') \\
c(w_1'r) &=& w_1'c(r) + c(w_1') \\
c(r^{-1}w_2') &=& r^{-1}c(w_2') + c(r^{-1})
\end{eqnarray*}
hold up to a uniformly bounded error. Hence, using the three equations above, we can compute:
\begin{eqnarray*}
c(w_1w_2) &=& c(w_1'w_2') \\
&=& w_1'c(w_2') + c(w_1') \\
&=& w_1' r c(r^{-1}w_2') - w_1' r c(r^{-1}) + c(w_1') \\
&=& w_1' r c(r^{-1}w_2') - w_1' r c(r^{-1}) + c(w_1'r') - w_1'c(r) \\
&=& w_1 c(w_2) + c(w_2)
\end{eqnarray*}
again up to uniformly bounded error. In the last step we used that $c$ is anti-symmetric. In fact, the construction
seems to fail at this point, if one does not assume $c_1$ and $c_2$ to be anti-symmetric, since we can not assure
that $c(r^{-1})+r^{-1}c(r)$ is uniformly bounded. This finishes the proof.
\end{proof}

\begin{remark}
Note that obviously, 
the class of $c$ in the first quasi-cohomology does depend heavily on $c_1$ and $c_2$ and 
not only on their classes in the first quasi-cohomology.
\end{remark}

\begin{theorem}
The class of groups $G$ which admit proper quasi-1-cocycles $c\colon G \to \ell^2G^{\oplus \infty}$ is closed under
subgroups and free products.
\end{theorem}
\begin{proof}
The assertion concerning subgroups is obvious since we can just restrict the cocycles and decompose the regular
representation according to the cosets. 

Let us now turn to the question about free products.
Given proper quasi-1-cocycles $c_1\colon G \to \ell^2G^{\oplus \infty}$
and $c_2\colon H \to \ell^2H^{\oplus \infty}$ we can regard both as taking values in $\ell^2(G \ast H)^{\oplus \infty}$ and can assume that they are anti-symmetric. Moreover, since the set of values of $c_1$ and $c_2$ in a bounded region is finite, we
can add a bounded $1$-cocycle and assume that the minimum of $g \mapsto \|c_1(g)\|_2$ and $h \mapsto \|c_2(h)\|$ is non-zero on $G \setminus \{e\}$ resp. $H \setminus \{e\}$.

We claim that the quasi-1-cocycle 
$(c_1 \ast 0) \oplus (0 \ast c_2)\colon G \ast H \to \ell^2(G \ast H)^{\oplus \infty}$ is proper.
Here, we are using the notation of Lemma \ref{freeprod}.
Let $w = g_1h_1g_2h_2 \dots g_n h_n$ be a reduced element in $G \ast H$. Clearly,
$$\|(c_1 \ast 0)(w)\|^2 = \sum_{i=1}^n \|c(g_i)\|^2.$$
Hence, using the properness of $c_1$, the set of $g_i$'s that can appear with a given bound on $\|(c_1 \ast 0)(w)\|$ is a finite subset of $G$. Moreover, we find an upper bound on $i$ since the minimum of $g \to \|c(g)\|$ is assumed to be non-zero. The same is true for the set of $h_i$'s. Hence, the set of elements $w \in G \ast H$, for which $\|(c_1 \ast 0) \oplus (0 \ast c_2)(w)\| $ is less than a constant is finite. This finishes the proof.
\end{proof}

\subsection{Orbit equivalence} \label{obeq}

In this section we study the stability of the class of groups which admit proper quasi-1-cocycles in a multiple of the
regular representation under orbit equivalence. For the notion of \emph{orbit equivalence}, which goes back to work of Dye (see \cite{dye1,dye2}), we refer to
the Gaboriau's nice survey \cite{MR2183295} and the references therein, see also \cite{MR1253544}.
Although it remains open, whether the class is closed under this relation,
we are able to prove that it is closed under a slightly more restricted relation, which we call $\ell^2$-orbit equivalence.
(To our knowledge, the idea of $\ell^2$-orbit equivalence goes back to unpublished work of R.\ Sauer.)
Unfortunately, we cannot say much more about the class of groups which are $\ell^2$-orbit equivalent to
Gromov hyperbolic groups. This is subject of future work. Literally everything extends to the suitable notions of
measure equivalence, but for sake of simplicity we restrict to orbit equivalence.

Let $G$ be a discrete countable group. Let $(X,\mu)$ be a standard probability space and $G \curvearrowright (X,\mu)$
be a measure preserving (m.p.) action by Borel automorphisms. 
We denote by $X \rtimes G$ the inverse semigroup of partial isomorphisms which are implemented by the action 
(not just the equivalence relation). 
Two partial isomorphisms $\phi,\psi$, which are induced by the action are said to be \emph{orthogonal} if $\dom(\phi)\cap \dom(\psi)=\varnothing$ and $\ran(\phi) \cap \ran(\psi) = \varnothing$. They are said to be \emph{disjoint} if they are disjoint
as subsets of the set of morphisms of the associated discrete measured groupoid. Clearly, orthogonal partial isomorphisms are disjoint. All equalities which concern subsets of a probability space or partial maps between probability spaces are supposed to hold almost everywhere, i.e.\ up to a set of measure zero, as usual.

Every partial isomorphism can be written as an infinite orthogonal sum as follows:
$$ \phi = \oplus_{i=1}^{\infty} \phi_{A_i} g_i,$$
for some Borel subsets $A_i$ and $g_i \in G$. 
The sub-inverse-semigroup of those for which there is a finite sum as above is denoted by $X \rtimes_{\rm fin}G$. If $G$ is finitely generated and $l\colon G \to \N$ is a word length function on $G$, there is yet another sub-inverse-semigroup, which is formed by those infinite sums, for which
$$\sum_{i=1}^{\infty} \mu(A_i) l(g_i)^2 < \infty.$$
We denote it by $X \rtimes_2 G$. Since $l(gh)^2 \leq (l(g) + l(h))^2 \leq 2 (l(g)^2 + l(h)^2)$, it is obvious that
$X \rtimes_2 G$ is closed under composition. Note also, that the summability does not depend on the set of generators
we choose to the define the length function.

\begin{definition} Let $(X,\mu)$ be a standard probability space and $G,H \curvearrowright (X,\mu)$ m.p.\ 
actions by Borel automorphisms. The data is said to induce an orbit-equivalence, if the orbits of the two actions 
agree up to measure zero. In this case, injective natural homomorphisms of inverse semigroups
$$\phi_1 \colon G \to X \rtimes H, \quad \mbox{and} \quad \phi_2\colon H \to X \rtimes G$$
are defined. We say that an orbit equivalence is an $\ell^2$-orbit-equivalence, if  the the image of $\phi_1$ (resp.\ $\phi_2$)
is contained in $X \rtimes_2 H$ (resp.\ $X \rtimes_2 G$).
\end{definition}

\begin{remark}
If the images are contained even in $X \rtimes_{\rm fin} G$ (resp.\ $X \rtimes_{\rm fin} H$, the one usually 
speaks about a uniform orbit equivalence.
Using Gromov's dynamical criterion, this also implies that $G$ is quasi-isometric to $H$. 
Hence, $\ell^2$-orbit equivalence is somehow 
half-way between quasi-isometry and usual orbit equivalence.
\end{remark}

\begin{definition}
Let $H_{\pi}$ be a unitary $G$-representation which carries a 
compatible normal action of $L^{\infty}(X)$. A 1-cocycle of $X \rtimes_2 G$ with 
values in $H_{\pi}$ is defined to be a map
$c \colon X \rtimes_2 G \to H_{\pi}$, such that
\begin{enumerate}
\item $c(\phi) \in \chi_{\ran(\phi)} H_\pi$,
\item $c$ is compatible with infinite orthogonal decompositions of the domain,
\item $c(\psi \phi) = \psi c(\phi) + c(\psi)$ if $\dom(\psi)=\ran(\phi)$.
\end{enumerate}
\end{definition}

A 1-cocycle is said to be inner, if $c(\phi)=(\phi - \chi_{\ran \phi}) \xi$, for some vector $\xi \in H_\pi$.
In analogy to the group case, we call a map $c\colon X \rtimes_2 G \to H_\pi$ satisfying 1. and 2. from above
a quasi-1-cocycle if $\| \psi c(\phi) - c(\psi \phi) + c(\psi)\|$ is uniformly bounded, for $\psi,\phi$ with $\dom(\psi)= \ran(\phi)$.

\begin{definition}
A quasi-1-cocycle $c\colon X \rtimes_2 G \to H_{\pi}$ is said to be \emph{proper}, if for every sequence of disjoint
partial isomorphisms $\phi_i \in X \rtimes G$ with $\liminf_{i \to \infty} \mu(\dom(\phi_i)) >0$, we have that
$\lim_{i \to \infty } \|c(\phi_i)\| = \infty$.
\end{definition}

\begin{lemma} \label{extension}
Let $G$ be a discrete countable group. Let $(X,\mu)$ be a standard probability space and $G \curvearrowright (X,\mu)$
be a m.p.\ action by Borel automorphisms. Let $c\colon G \to H_{\pi}$ be a quasi-1-cocycle. There is a natural
extension of the (quasi-)1-cocycle $c$ to a (quasi-)1-cocycle 
$\tilde{c}\colon X \rtimes_{2} G \to L^2(X,\mu) \otimes_2 H_{\pi}$ where $G$ acts diagonally. 
\end{lemma}

\begin{proof}
We define $$\tilde{c}\left(\oplus_{i=1}^\infty \chi_{A_i} g_i \right)= \sum_{i=1}^\infty \chi_{A_i} \otimes c(g_i).$$
Since $\|c(g)\|\leq C\cdot l(g)$ for some constant $C>0$, the right hand side is well-defined in $L^2(X,\mu) \otimes_2 H_{\pi}$. It can be easily checked that all relations are satisfied.
\end{proof}

\begin{lemma} \label{proper}
If the quasi-1-cocycle $c \colon G \to H_\pi$ is proper, then so is the quasi-1-cocycle 
$\tilde{c}\colon X \rtimes_2 G \to L^2(X,\mu) \otimes_2 H_{\pi}$, which we obtain from the construction 
in Lemma \ref{extension}.
\end{lemma}
\begin{proof}
Let $\phi_i$ be a sequence of disjoint partial isomorphisms with $\liminf_{i \to \infty} \mu(\phi_i) \geq \varepsilon > 0$. 
In order to derive a contradiction, we can assume that $\|\tilde{c}(\phi_i)\| < C$ for some constant $C$ and all $i \in \N$.

Hence, for every $i \in \N$, at least half of $\phi_i$ is supported at group elements $g \in G$ 
with $\|c(g)\| \leq 2C/\varepsilon$. Indeed, if otherwise, then 
$\|\tilde{c}(\phi_i)\| \geq 2^{-1/2} \varepsilon \cdot 2C/\varepsilon >C$. Since this holds for all $i \in \N$, and the support of 
$g \in G$ has measure $1$, the set of $g \in G$, with $\|c(g)\| \leq 2C/\varepsilon$ has to be infinite. This is a contradiction, since we assume $c$ to be proper.
\end{proof}

\begin{theorem}\label{main3}
The class of groups which admit a proper quasi-1-cocycle with values in an infinite sum of the regular representation
is closed under $\ell^2$-orbit equivalence.
\end{theorem}
\begin{proof} 
Let $G$ and $H$ be $\ell^2$-orbit equivalent groups. We show that, if $G$ admits a proper quasi-1-cocycle with values in $\ell^2 G^{\oplus \infty}$,
then $H$ admits a proper quasi-1-cocycle with values in $\ell^2H^{\oplus \infty}$.

Let $c\colon G \to \ell^2 G^{\oplus \infty}$ be a proper quasi-1-cocycle. 
Lemma \ref{extension} says that we can extend $c$ to a quasi-1-cocycle which is defined on $X \rtimes_2 G$ and takes
values in $L^2(X,\mu) \otimes_2 \ell^2 G^{\oplus \infty}$. Note that the homomorphism 
$\phi\colon H \to X \rtimes_2 G$ is compatible, in the sense that the obvious actions are intertwined with the natural isomorphism 
$$L^2(X,\mu) \otimes_2 \ell^2 G^{\oplus \infty} \cong L^2(X,\mu) \otimes_2 \ell^2 H^{\oplus \infty}.$$

The proof is finished by noting that the restriction $\tilde{c}|_H\colon H \to L^2(X,\mu) \otimes_2 \ell^2 H^{\oplus \infty}$
is proper by Lemma \ref{proper}, and that $L^2(X,\mu) \otimes_2 \ell^2 H^{\oplus \infty} \cong \ell^2 H^{\oplus \infty}$
as unitary $H$-representations.
\end{proof}

\begin{remark}
A similar proof applies to the case where $H$ merely embeds into $X \rtimes_2 G$ (i.e. is not necessarily part of
an orbit equivalence). This type of subobject is called \emph{random subgroup} in \cite{MR2275641}. However, note that 
an arbitrary random subgroup does not necessarily satisfy the $\ell^2$-condition we impose. We could speak of
$\ell^2$-random subgroups in case it satisfies the $\ell^2$-condition.
\end{remark}

Note that the following strengthening of Theorem \ref{main2} is also an immediate consequence of Theorem \ref{main3}.
\begin{corollary} \label{coromain2}
Let $G$ be a group which is boundedly generated by a finite set of amenable subgroups and let $H$ be
a Gromov hyperbolic group. Let $H \curvearrowright (X,\mu)$ be a m.p.\ action by Borel automorphisms on
a standard probability space. Any homomorphism of inverse semi-groups $\phi\colon G \to X \rtimes_2 H$
has amenable image. Moreover, if $G$ has property (T) of Kazhdan, then the image of $\phi$ has finite measure.
\end{corollary}

\begin{remark}
In terms of randomorphisms, the last corollary says that any $\ell^2$-randomorphism from $G$ to $H$ has 
finite image; in the sense that the relevant $G$-invariant measure $\mu$ on the polish space $[G,H]_\bullet$ 
(see \cite{MR2275641} for details) is a.e.\ supported on maps with finite image such that 
$$\int_{[G,H]_\bullet} \# \phi(G)\ d\mu(\phi) < \infty.$$
\end{remark}

\medskip

\begin{center}
{\bf Acknowledgments}
\end{center}

Thanks go to Theo B\"uhler, Clara L\"oh and Thomas Schick for interesting discussions about bounded cohomology.
I am grateful to Roman Sauer for sharing his ideas about $\ell^2$-orbit equivalence.
I also want to thank the unknown referee for several useful remarks that improved the exposition.

\end{document}